\documentclass{article}
\usepackage{latexsym, amssymb}
\newtheorem{thm}{Theorem}

\newtheorem{cor}{Corollary}
\newtheorem{prop}{Proposition}
\newtheorem{lem}{Lemma}
\newtheorem{defn}{Definition}
\newtheorem{quest}{Question}

\begin{document}

\begin{center}
{\bf FINITE TYPE LINK HOMOTOPY INVARIANTS II:  Milnor's $\bar{\mu}$-invariants}\\
\vspace{.2in}
{\footnotesize BLAKE MELLOR}\\
{\footnotesize University of California, Berkeley}\\
{\footnotesize Department of Mathematics}\\
{\footnotesize Berkeley, CA  94720}\\
{\footnotesize\it  mellor@math.berkeley.edu}\\
\vspace{.2in}
{\footnotesize\sc December 1998}\\
\vspace{.2in}
{\footnotesize\it This is a preprint - all comments are welcome!}\\
\vspace{.2in}
{\footnotesize ABSTRACT}\\
{\ }\\
\parbox{4.5in}{\footnotesize \ \ \ \ \ We define a notion of finite type invariants for links 
with a fixed linking matrix.  We show that Milnor's link homotopy invariant $\bar{\mu}(ijk)$ 
is a finite type invariant, of type 1, in this sense.  We also generalize this approach to 
Milnor's higher order $\bar{\mu}$ invariants and show that they are also, in a sense, of 
finite type.  Finally, we compare our approach to another approach for defining finite type 
invariants within linking classes.}\\
\vspace{1in}
\end{center}

\input{psfig.sty}

\tableofcontents

\section{Introduction} \label{S:intro}

In the usual theory of finite type invariants (see \cite{bn1} for a thorough introduction), 
we use crossing changes to move between isotopy classes of links.  These classes are components 
of the space of all embeddings $\sqcup S^1 \hookrightarrow S^3$, which is contained within 
the space of all immersions $\sqcup S^1 \rightarrow S^3$.  The boundaries of the isotopy 
classes (in the space of immersions) are links with double points, and we extend link invariants 
linearly to these boundaries by the formula:
$$\psfig{file=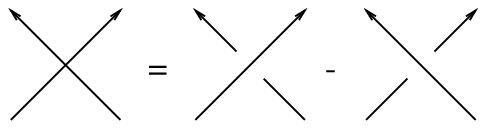}$$
A {\it finite type} invariant of {\it type m} is then an extended invariant which vanishes 
on links with more than $m$ double points.  However, there are invariants which are not 
well-defined on these differences, and so cannot be analyzed within the usual framework of finite 
type invariants.  For example, Milnor's link-homotopy invariant $\bar{\mu}(ijk)$ is only 
well-defined modulo the linking numbers of the three components; since crossing changes alter 
these linking numbers, the invariant cannot be meaningfully extended to the boundaries 
between link types.  In fact, the author has shown in \cite{me} that, up to link-homotopy, 
the only finite type invariants in the usual sense are the linking numbers.  The goal of this 
paper is to define a notion of finite type invariant {\it within} a class of links with the 
same linking matrix, where $\bar{\mu}(ijk)$ is well-defined.  We will call the class of all 
links with a given linking matrix a {\it linking class}.  This work was inspired by recent work 
of Cochran and Melvin on finite type invariants of 3-manifolds (see \cite{cm}), where they 
used similar ideas to define a notion of finite type invariants for all 3-manifolds, rather 
than just homology 3-spheres.  

In section~\ref{S:borrmu} we will introduce the notion of a {\it Borromean clasp} and define a 
theory of finite type invariants based on adding and removing these clasps (an operation which 
preserves the linking class).  We show that $\bar{\mu}(ijk)$ is of type 1 in this theory.  In 
section~\ref{S:higher} we generalize this idea to the higher $\bar{\mu}$-invariants and show 
that they are also, in a sense, type 1 invariants.

In section~\ref{S:claspclass} we describe the equivalence classes of the clasping operations 
defined in the previous sections.  In particular, for the Borromean clasps these equivalence 
classes are the linking classes.

In section~\ref{S:double-crossing} we consider another theory of finite type invariants for 
linking classes, based on changing pairs of oppositely signed crossings.  This theory is also 
being studied by Eli Appleboim and Dror Bar-Natan (see [A-BN]).  We show that 
$\bar{\mu}(ijk)$ is also of finite type in this theory, though of type 2 rather than type 1.  
Finally, in section~\ref{S:compare}, we compare the two theories we have developed for linking 
classes, and find that the Borromean clasp theory may be stronger (though not necessarily more 
useful - see the questions in section~\ref{S:questions}).

\section{Borromean clasps and $\bar{\mu}(ijk)$} \label{S:borrmu}

\subsection{Borromean clasps} \label{SS:borr}

We will begin by considering links up to link homotopy.  
The basic idea is to look at a different notion of ``crossing change''.  The usual notion of 
crossing changes can be thought of as removing a ``clasp'' between two strands, as shown 
below:
$$\psfig{file=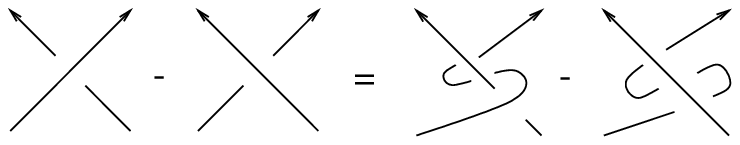}$$
We want to find an analogous clasp for linking classes up to link homotopy.  The obvious 
choice is the Borromean rings, which has trivial linking numbers but is well-known to be 
homotopically non-trivial.  So we will look at the following ``crossing change'':
$$\psfig{file=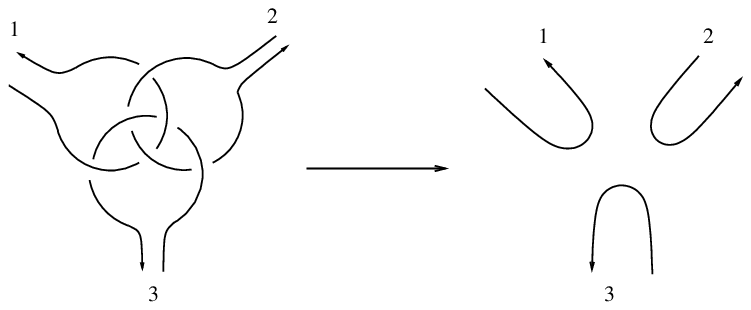}$$
Clearly, this operation preserves the linking class of the link.  It will be useful to 
``straighten out'' this clasp and look at it as an operation on string links.  In this 
case, the operation is:
$$\psfig{file=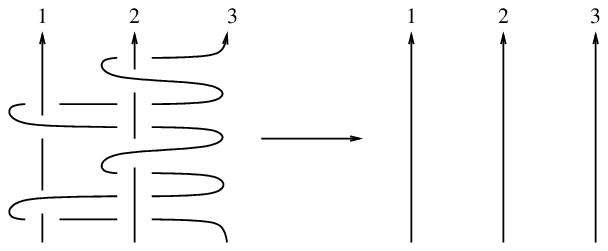}$$
From \cite{hl}, we know that the group of string links up to link homotopy (under the 
operation of concatenation), is generated by $\{x_{ij}|i<j\}$:
$$\psfig{file=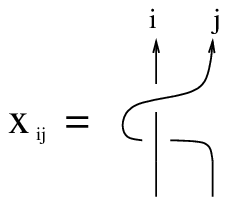}$$
Then, in terms of these generators, the Borromean clasp is the commutator $[x_{13},x_{23}]$, 
and removing the clasp corresponds to removing this commutator from the word in the $x_{ij}$'s 
corresponding to the string link.  This suggests that perhaps we should generalize our 
idea of the ``Borromean clasp'' to include all commutators $[x_{ij}^{\pm 1},x_{kl}^{\pm 1}]$.  
These commutators are trivial if $i,j,k,l$ are all distinct, so it is not hard to see that it 
suffices to consider the commutators $[x_{ik}^{\pm 1},x_{jk}^{\pm 1}]^{\pm 1}$, where $i<j<k$.  
Hence, we are considering all the operations:
$$\psfig{file=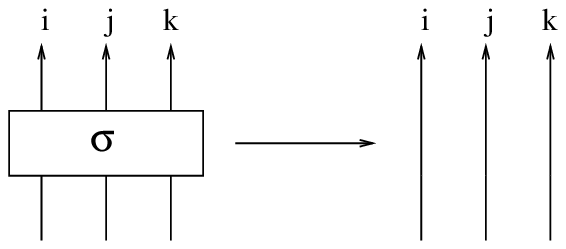}$$
Where $\sigma = [x_{ik}^{\pm 1},x_{jk}^{\pm 1}]^{\pm 1}$.

Clearly, these all fix the linking class, since the linking number of components $i$ and $j$ 
is simply the multiplicity, with sign, of the generator $x_{ij}$, which is not changed by 
adding or removing a commutator.  Now we can define the notions of a singular link and a 
finite type invariant as in the usual theory.

\begin{defn} \label{D:borr-singular}
A {\bf singular link} of degree m (in the Borromean clasp theory) is a link with m triple 
points, each labeled by a commutator $\sigma = [x_{ik}^{\pm 1},x_{jk}^{\pm 1}]^{\pm 1}$:
$$\psfig{file=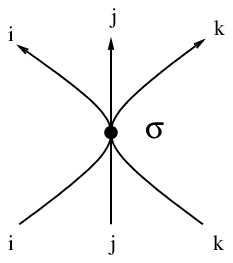}$$
\end{defn}

We extend any link-homotopy invariant which is well-defined within each linking class 
to singular links by the relation:
$$\psfig{file=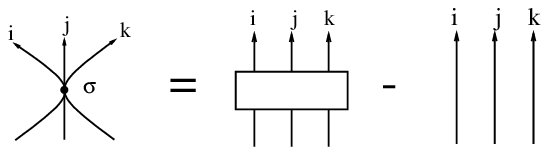}$$

\begin{defn} \label{D:borr-type}
An invariant V defined within a linking class is of {\bf type m} if it is trivial on all 
singular links (within that linking class) of degree $\geq m+1$.  V is of {\bf finite type} 
if it is of type m for some finite m.
\end{defn}

Now we want to show that $\bar{\mu}(ijk)$ is of finite type in this sense.  First, we will 
review the definition of Milnor's $\bar{\mu}(ijk)$ invariants.

\subsection{$\bar{\mu}(ijk)$} \label{SS:muijk}

We recall how to compute $\bar{\mu}(ijk)$ from \cite{mi}.  Given a link $L$, its link group 
$\pi_1(S^3 - L)$ has a Wirtinger presentation, generated by the arcs of the link diagram.  We 
also have a presentation of the link group modulo the $q$th subgroup in its lower central series 
(see \cite{mi}):
$$\pi_1(S^3-L)/(\pi_1)_q(S^3-L) = <m_i\ |\ m_il_im_i^{-1}l_i^{-1} = 1, A_q>$$
where the generators are the meridians $m_i$ of the components of the link, the $l_i$ denote 
the longitudes of the components of the link, and $A_q$ denotes the $qth$ subgroup in the lower 
central series of the free group on $\{m_i\}$.  So each longitude (and the generators of the 
Wirtinger presentation) can be written in $\pi_1/(\pi_1)_q$ as a word in the $m_i$'s.  We 
look at the Magnus expansion of the longitudes, which means replacing $m_i$ with $1+K_i$ and 
$m_i^{-1}$ with $1-K_i+K_i^2-...$.  We define $\mu(ijk)$ as the coefficient of $K_iK_j$ in the 
Magnus expansion of $l_k$.  In general, this is not well-defined for links.  Then 
$\bar{\mu}(ijk)$ is $\mu(ijk)$ modulo $\Delta = gcd\{linking\ numbers\ for\ components\ 
i, j, k\}$.  This is now a well-defined invariant of links up to concordance, 
as long as $q>2$ (it is otherwise independent of $q$).  If the indices $i,j,k$ are all 
distinct, it is in fact an invariant of link homotopy (see \cite{mi}).
So, within each linking class, we can extend $\bar{\mu}(ijk)$ to singular links.  Notice 
that $\bar{\mu}(ijk)$ will be trivial on any singular link with singularities involving a 
component other than $i,j,k$.

Note that $\mu(ijk)$, while not well-defined for links, {\it is} well-defined 
for string links (also see \cite{hl}).  Since the components of the string link have a 
natural beginning and ending, there are no choices made in producing the word in the 
meridians associated with each longitude, and it is easy to check that this word is 
invariant under the Reidemeister moves and string link homotopy.  In fact, modding out 
by $\Delta$ exactly compensates for the effect of closing a string link to get a link, 
as the following lemma shows:

\begin{lem} \label{L:closing}
If L is a link and $\sigma$ is a string link such that L = $\hat{\sigma}$, then 
$\bar{\mu}(ijk)(L) = \mu(ijk)(\sigma)\ mod\ \Delta$.
\end{lem}
{\sc Proof:}  Consider a component $l_r$ of $\sigma$, and let $\alpha_r$ and $\gamma_r$ denote the 
meridians of $l_r$ at its beginning and end, respectively.  By abuse of notation, we will also 
let $l_r$ denote the word in the $\alpha_r$'s representing the longitude of $l_r$.  Then it is 
clear that $\gamma_r = l_r^{-1}\alpha_r l_r$.  The effect of closing the string link $\sigma$ 
is to identify $\gamma_r$ with $\alpha_r$.  How does this identification change $\mu(ijk)$?  
The contribution to $\mu(ijk)$ of each appearance of $\gamma_r$ in $l_k$ is almost the same 
as that of an appearance of $\alpha_r$, with the difference being the coefficient of $K_iK_j$ 
in the Magnus expansion of $\gamma_r$, which is ($mult_i(w)$ is the multiplicity, with sign 
of $m_i$ in the word $w$):
$$mult_i(l_r^{-1})mult_j(\alpha_r) + mult_i(\alpha_r)mult_j(l_r)$$
$$ = -lk(l_i, l_r)\delta(j,r) + \delta(i,r)lk(l_j, l_r)$$
which is $\pm lk(l_i,l_j)$ if $r = i,j$, and 0 otherwise.  Modulo $\Delta$, this difference 
disappears, so the effect of $\gamma_r$ is exactly the same as that of $\alpha_r$, and we 
have computed $\bar{\mu}(ijk)(\hat{\sigma}) = \bar{\mu}(ijk)(L)$.  $\Box$

Because of this lemma, we can do some of our analysis for string links instead of links.  In 
particular, since string links have a group structure, we can look at how $\mu(ijk)$ behaves 
under multiplication.  The proof of the following lemma is immediate:

\begin{lem} \label{L:additivity}
If $\sigma$, $\sigma_1$ and $\sigma_2$ are string links such that $\sigma = \sigma_1\sigma_2$, 
then $\mu(ijk)(\sigma) = \mu(ijk)(\sigma_1) + \mu(ijk)(\sigma_2) + \mu(ik)(\sigma_1)\cdot
\mu(jk)(\sigma_2)$ (where $\mu(ij)$ is just the linking number of the $i$th and $j$th 
components, see \cite{mi}).  In particular, if $\sigma_1$ or $\sigma_2$ are algebraically 
unlinked, then $\mu(ijk)$ acts additively.
\end{lem}

\subsection{$\bar{\mu}(ijk)$ is finite type} \label{SS:ijk-borr-finite}

Now we can use our understanding of $\bar{\mu}(ijk)$ to prove that it is a finite type 
invariant (in the Borromean clasp theory).

\begin{thm} \label{T:ijk-type1}
$\bar{\mu}(ijk)$ is of type 1 (in the Borromean clasp theory).
\end{thm}
{\sc Proof:}  We need to know what the contribution of a commutator 
$[x_{ik}^{\pm 1}, x_{jk}^{\pm 1}]^{\pm 1}$ is to $\mu(ijk)$ on the level of string links.  
Say that $\gamma$ is such a commutator, and $\sigma = \delta\gamma$ is a string link.  
Then, by Lemma 2, $\mu(ijk)(\sigma) = \mu(ijk)(\delta) + \mu(ijk)(\gamma) + 
\mu(ij)(\delta)\mu(jk)(\gamma)$.  Since $\mu(jk)(\gamma) = 0$, we conclude that the 
contribution of $\gamma$ to $\mu(ijk)(\sigma)$ is exactly $\mu(ijk)(\gamma)$.  In particular, 
the contribution depends only on $\gamma$, and is independent of the rest of the link.  

So we need to look at the Magnus expansion of $[x_{ik}^{\pm 1}, x_{jk}^{\pm 1}]^{\pm 1}$ 
(along component $k$), and find the coefficient of $K_iK_j$.  First we consider the commutator 
$[x_{ik},x_{jk}] = x_{ik}x_{jk}x_{ik}^{-1}x_{jk}^{-1}$.  This contributes a word 
$m_im_jm_i^{-1}m_j^{-1}$ to the word for the $k$th component.  The Magnus expansion is 
thus $(1+K_i)(1+K_j)(1-K_i+...)(1-K_j+...) = 1 + K_iK_j - K_iK_j + K_iK_j + ... = 1 + 
K_iK_j + ...$, so $\mu(ijk)([x_{ij},x_{jk}]) = 1$.  Similarly, for the other commutators, 
$\mu(ijk) = \pm 1$.

Now, say that $L$ is a link with two Borromean singularities $a$ and $b$.  $L$ is the 
closure of a singular string link $\sigma$.  If either of the singularities $a$ or $b$ involve 
a component other than $i,j,k$, then $\mu(ijk)(\sigma) = 0$.  Otherwise, $\sigma$ can be 
written as a linear combination of four string links:  $\sigma = \sigma' - \sigma'_a - 
\sigma'_b + \sigma'_{ab}$, where the subscripts indicate which of the Borromean clasps have 
been unclasped.  Then $\mu(ijk)(\sigma) = (-1)^a - (-1)^a = 0$ (by abuse of notation, $(-1)^a$ 
is $\mu(ijk)$ of the commutator corresponding to the singularity $a$).  Then by Lemma 1, 
$\bar{\mu}(ijk)(L)$ is also 0, which proves that $\bar{\mu}(ijk)$ is a finite type invariant 
of type 1 (in the Borromean clasp theory).  $\Box$

\section{Higher Milnor invariants} \label{S:higher}

In this section we will generalize the results for $\bar{\mu}(ijk)$ to higher $\bar{\mu}$-
invariants.  The higher invariants are defined in the same way as $\bar{\mu}(ijk)$.  First 
we compute $\mu(i_1...i_n,j)$, which is the coefficient of $K_{i_1}...K_{i_n}$ in the 
Magnus expansion of the word for the $j$th longitude in $\pi_1/(\pi_1)_q,\ q>n$, and then we 
consider it modulo $\Delta$, which is the greatest common divisor of all $\mu$-invariants 
whose indices are a cyclic permutation of a proper subsequence of $(i_1...i_nj)$.  
If the indices are all distinct, this is a well-defined link-homotopy invariant.

Now we can easily prove the analogues to Lemmas~\ref{L:closing} and \ref{L:additivity}:

\begin{lem} \label{L:higherclosing}
If L is a link and $\sigma$ is a string link such that L = $\hat{\sigma}$, then 
$\bar{\mu}(i_1...i_n,j)(L) = \mu(i_1...i_n,j)(\sigma)\ mod\ \Delta$.
\end{lem}
{\sc Proof:}  Consider a component $l_r$ of $\sigma$, and let $\alpha_r$ and $\gamma_r$ denote the 
meridians of $l_r$ at its beginning and end, respectively.  By abuse of notation, we will also 
let $l_r$ denote the word in the $\alpha_r$'s representing the longitude of $l_r$.  Then it is 
clear that $\gamma_r = l_r^{-1}\alpha_r l_r$.  The effect of closing the string link $\sigma$ 
is to identify $\gamma_r$ with $\alpha_r$.  How does this identification change $\mu(i_1...i_n,j)$?  
The contribution to $\mu(i_1...i_n,j)$ of each appearance of $\gamma_r$ in $l_k$ is almost the same 
as that of an appearance of $\alpha_r$, with the difference being the coefficient of $K_{i_1}...K_{i_n}$ 
in the Magnus expansion of $\gamma_r$, which is ($mult_I(w)$ is the coefficient of $K_I$ in the Magnus 
expansion of $w$):
$$\sum_{k=1}^n{mult_{i_1...i_{k-1}}(l_r^{-1})mult_{i_k}(\alpha_r)mult_{i_{k+1}...i_n}(l_r)}$$
$$= \sum_{k=1}^n{-mult_{i_{k-1}...i_1}(l_r)mult_{i_k}(\alpha_r)mult_{i_{k+1}...i_n}(l_r)}$$
$$= \sum_{k=1}^n{-\mu(i_{k-1}...i_1,r)\mu(i_{k+1}...i_n,r)\delta(i_k,r)}$$
$$= \sum_{k=1}^n{(-1)^k\mu(i_1...i_{k-1},r)\mu(i_{k+1}...i_n,r)\delta(i_k,r)}$$
The last equality is a result of Milnor, see \cite{mi}.  
This sum is 0 unless $r=i_k$ for some $k$, in which case the sum is equal to the $k$th term, i.e. 
$(-1)^k\mu(i_1...i_{k-1},i_k)\mu(i_{k+1}...i_n,i_k)$.  But this is trivial modulo $\Delta$, so $\gamma_r$ 
has the same effect as $\alpha_r$, and we conclude that $\bar{\mu}(i_1...i_n,j)(L) = 
\mu(i_1...i_n,j)(\sigma)$.  $\Box$

\begin{lem} \label{L:higheradditivity}
If $\sigma$, $\sigma_1$ and $\sigma_2$ are string links such that $\sigma = \sigma_1\sigma_2$, 
then:
$$\mu(i_1...i_n,j)(\sigma) = \mu(i_1...i_n,j)(\sigma_1) + \mu(i_1...i_n,j)(\sigma_2)$$
$$+ \sum_{k=1}^{n-1}{\mu(i_1...i_k,j)(\sigma_1)\cdot\mu(i_{k+1}...i_n,j)(\sigma_2)}$$
\end{lem}
{\sc Proof:}  The proof is trivial.  $\Box$

Now we need to determine the appropriate ``clasps'' for the higher order Milnor invariants.  Since 
these invariants are only defined modulo the lower order invariants, we need to find operations 
which preserve the lower order invariants.  Since the Milnor invariants are, in a sense, higher 
order linking numbers (see \cite{co} for a geometric approach to the Milnor invariants which 
makes this precise), we will call a class of links which share all Milnor invariants of length 
$n$ or less an {\it n-linking class}.  In this terminology, our usual linking classes are 2-linking 
classes.  For 2-linking classes, our clasps were commutators, elements of the second group in the 
lower central series.  For $n$-linking classes, we will look at the $n$th group of the lower 
central series.  This is a sensible approach, since Cochran shows that Milnor's invariants 
are really measuring how deep each longitude lies in the lower central series of the link group 
(see \cite{co}).

To make this precise, let $H(k)$ denote the group of string links with $k$ components (the 
group operation is just concatenation - put one string link ``on top'' of the other).  The lower 
central series is defined inductively by $H_1(k) = H(k)$ and $H_n(k) = [H(k),H_{n-1}(k)]$.  
Habegger and Lin have shown (see \cite{hl}) that $H(k)$ is nilpotent of order $k-1$; in other 
words, $H_k(k) = 1$.  As before, $H(k)$ is generated by $\left(\matrix{k \cr 2}\right)$ generators 
$x_{ij},\ i<j$.  We now define the following elements of the lower central series:

\begin{defn} \label{D:simplecommutator}
Given n components (WLOG, numbered 1 to n), we define a {\bf simple 1-commutator} as a generator 
$x_{in}^{\pm 1}$.  We then inductively define a {\bf simple k-commutator} as an element of 
$H_k(n)$ of the form $[x_{in}^{\pm 1}, A]$ or $[A,x_{in}^{\pm 1}]$, where A is a simple (k-1)-
commutator.
\end{defn}

Notice that a simple $k$-commutator is a Brunnian link with $k+1$ components (i.e. removing any of 
the components trivializes the link).  Therefore, any $\mu$-invariant of length $k$ or less is 
trivial on any $k$-commutator, since it can only ``see'' at most $k$ of the components.  
Our operations in an $n$-linking class will now consist of removing simple $n$-commutators:
$$\psfig{file=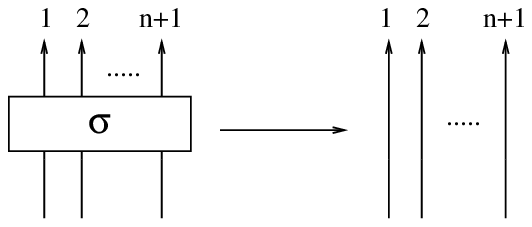}$$
Now we can generalize our earlier definitions of singular links and finite type invariants.

\begin{defn} \label{D:n-singular}
An {\bf n-singular link} of degree m (in an n-linking class) is a link with m (n+1)-tuple 
points, each labeled with a simple n-commutator:
$$\psfig{file=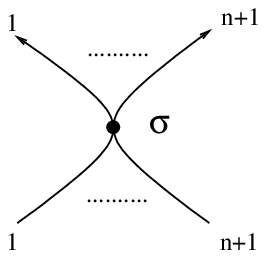}$$
\end{defn}

We extend any link-homotopy invariant which is well-defined withing each $n$-linking class 
to $n$-singular links by the relation:
$$\psfig{file=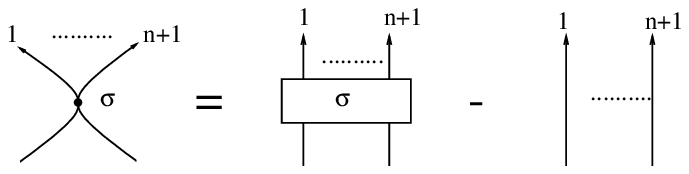}$$

\begin{defn} \label{D:n-type}
An invariant V defined within an n-linking class is of {\bf n-type m} if it is trivial on all 
n-singular links (within that n-linking class) of degree $\geq m+1$.  V is of {\bf finite n-type} 
if it is of n-type m for some finite m.
\end{defn}

Now we want to show that the Milnor invariants of length $n+1$ are all of finite $n$-type.  The 
first step is to look at the values of these invariants on simple $n$-commutators, on the level 
of string links.

\begin{lem} \label{L:mu(commutator)}
Let $\sigma$ be a simple n-commutator on components $1,...,n+1$.  Then 
$\mu(1,...,n+1)(\sigma) = 0,\ 1,\ or\ -1$.
\end{lem}
{\sc Proof:}  The proof is by induction on $n$.  The lemma is true for $n = 2$ by the proof of 
Theorem 1.  So assume it is true for $n-1$.  First note that if any $x_{i(n+1)}$ does not appear in 
$\sigma$, then $\sigma$ only involves the $n$ components $1,...,i-1,i+1,...,n+1$, which means that 
$\sigma \in H_n(n) = 1$, so $\sigma$ is trivial, and $\mu(1,...,n+1)(\sigma)=0$.  So we may assume 
each $x_{i(n+1)}$ appears exactly once in $\sigma$.  Then $\sigma = [x_{i(n+1)}^{\pm 1},A]$ or 
$[A,x_{i(n+1)}^{\pm 1}]$, where $A$ is a simple $(n-1)$-commutator on $\{x_{1(n+1)},...,x_{(i-1)(n+1)},
x_{(i+1)(n+1)},...,x_{n(n+1)}\}$.  

We will assume $\sigma = [x_{i(n+1)},A]$ (the other cases are similar).  Now we apply 
Lemma~\ref{L:higheradditivity} to compute $\mu(1...n+1)(\sigma)$.  Notice that the only non-trivial 
$\mu$-invariants of $x_{i(n+1)}$ are of length 2, and the only non-trivial $\mu$-invariants of $A$ 
are of length $n$ (since $A$ is a simple $n-1$-commutator).  So we drop all other terms without further 
comment.
$$\mu(1...n+1)(\sigma) = \mu(1...n+1)(x_{i(n+1)}Ax_{i(n+1)}^{-1}A^{-1})$$
$$ = \mu(1...n+1)(Ax_{i(n+1)}^{-1}A^{-1}) + \mu(1,n+1)(x_{i(n+1)})\mu(2...n+1)(Ax_{i(n+1)}^{-1}A^{-1})$$
$$ = \mu(1...n-1,n+1)(A)\mu(n,n+1)(x_{i(n+1)}^{-1}A^{-1}) + \mu(1...n+1)(x_{i(n+1)}^{-1}A^{-1})$$
$$ + \mu(1,n+1)(x_{i(n+1)})\mu(2...n+1)(A) + \mu(1,n+1)(x_{i(n+1)})\mu(2...n+1)(x_{i(n+1)}^{-1}A^{-1})$$
$$ = \mu(1...n-1,n+1)(A)\mu(n,n+1)(x_{i(n+1)}^{-1}) + \mu(1,n+1)(x_{i(n+1)}^{-1})\mu(2...n+1)(A^{-1})$$
$$ + \mu(1,n+1)(x_{i(n+1)})\mu(2...n+1)(A) + \mu(1,n+1)(x_{i(n+1)})\mu(2...n+1)(A^{-1})$$
$$ = \mu(1...n-1,n+1)(A)\mu(n,n+1)(x_{i(n+1)}^{-1}) + \mu(1,n+1)(x_{i(n+1)})\mu(2...n+1)(A)$$
$$ = \left\{\matrix{-\mu(1...n-1,n+1)(A)\ if\ i=n \cr 
			\mu(2...n+1)(A)\ if\ i=1 \cr 
			0\ otherwise}\right.$$
By induction, this is 0 or $\pm 1$, which completes the proof of the lemma.  $\Box$

\begin{thm} \label{T:higher-type1}
$\bar{\mu}(i_1...i_n,j)$ is of $n$-type 1.
\end{thm}
{\sc Proof:}  The proof is almost identical to the proof of Theorem 1, using the analogous lemmas for 
the higher-order invariants.  $\Box$

\section{Equivalence Classes of the Clasping Operations} \label{S:claspclass}

A theory of finite type isn't very useful unless the operation it is based on actually changes links.  
After all, if the operation is a Reidemeister move, then any invariant is of type 0!  In this section, 
we will look at the equivalence relation defined by the clasping operations described in the earlier 
sections of the paper.  Two links in the same $n$-linking class will be considered {\it equivalent} 
if one can be transformed into the other by adding and removing $n$-clasps - i.e. the clasps 
corresponding to simple $n$-commutators.  Ideally, the equivalence classes would be equal to the 
$n$-linking classes.  We show below that this is true for $n=2$.  For $n>2$ we have possibly 
smaller equivalence classes (in section~\ref{S:questions} we ask whether they are in fact smaller).

First, we state a useful lemma (see \cite{kms} for a proof):

\begin{lem} \label{L:Witt-Hall}
(Witt-Hall identities)  Let G be a group and let k, l, m be positive integers.  Say that $x\in G_k$, 
$y\in G_l$ and $z\in G_m$.  Then we have the following properties:
\begin{enumerate}
	\item  $[G_k,G_l] \subset G_{k+l}$, or $xy \equiv yx\ mod\ G_{k+l}$
	\item  $[x,zy] = [x,z][x,y][[y,x],z]$
	\item  $[x,zy] = [y,z][[z,y],x][x,z]$
	\item  $[x,[y,z]][y,[z,x]][z,[x,y]] \equiv 1\ mod\ G_{k+l+m+1}$
	\item  If $g \equiv g'\ mod\ G_k$ then $[g,y] \equiv [g',y]\ mod\ G_{k+l}$ and $[y,g] \equiv 
[y,g']\ mod\ G_{k+l}$.
\end{enumerate}
\end{lem}

\begin{thm} \label{T:claspclass}
Let L and $L'$ be two links with k components.  L can be transformed to $L'$ (up to link homotopy) 
by adding or removing simple 
n-commutators $\Leftrightarrow$ there exist string links $\sigma$ and $\sigma'$ such that 
$\hat{\sigma}$ = L, $\hat{\sigma}'$ = $L'$, and $\sigma \equiv \sigma'$ modulo $H_n(k)$.
\end{thm}
{\sc Proof:}  $(\Rightarrow)$  Let $\sigma$ be a string link with closure $L$.  Then each move on $L$ 
corresponds to inserting a simple $n$-commutator (or its inverse) into $\sigma$, and each simple 
$n$-commutator is an element of $H_n(k)$.  The result is a string link $\sigma'$ which is congruent 
to $\sigma$ modulo $H_n(k)$, and whose closure is $L'$.

$(\Leftarrow)$  $\sigma = \gamma\sigma'$, where $\gamma \in H_n(k)$.  By the Witt-Hall identities, 
$\gamma$ can be written as a product of simple $n$-commutators, modulo $H_{n+1}(k)$.  We can continue 
this process at each level; since $H(k)$ is nilpotent, it will terminate.  So we can write $\gamma$ 
as a finite product of simple $l$-commutators, where $n \leq l \leq k-1$.  In each of these commutators, 
there is an innermost simple $n$-commutator which can be removed (in each of its appearances), 
causing the larger commutator to disappear.  Since there are a finite number of commutators, each 
of finite length, we conclude that removing a finite number of simple $n$-commutators (each 
corresponding to an unclasping operation) will trivialize $\gamma$, and transform $\sigma$ 
to $\sigma'$.  Hence these moves will also transform $L$ to $L'$.  $\Box$

\begin{cor} \label{C:borr-class}
L can be transformed to $L'$ by adding or removing Borromean clasps (simple 2-commutators) 
$\Leftrightarrow$ L and $L'$ are in the same linking class.
\end{cor}
{\sc Proof:}  Linking numbers completely determine string links up to homotopy modulo $H_2(k)$ 
(see \cite{hl}).  $\Box$

\section{Double Crossing Changes} \label{S:double-crossing}

There is another, perhaps more obvious, operation on links which fixes the linking class; namely, 
we change pairs of crossings with opposite sign, as shown below:
$$\psfig{file=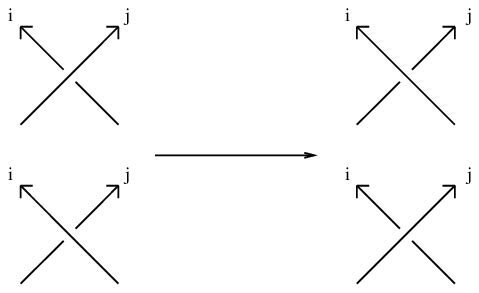}$$
In this section we will look at the finite type theory generated by this operation.  This has also 
been studied by other authors, in particular by Appleboim and Bar-Natan (see \cite{bna}).  We will 
show that $\bar{\mu}(ijk)$ is also of finite type in this theory, though here it is of type 2.  In 
the next section, we will compare this theory with the theory generated by Borromean clasps.

\subsection{Definitions and Equivalence classes of the Double Crossing Changes}

We define a notion of finite type invariants in the context of double crossing changes.  We begin 
by defining a notion of singular link (see also \cite{bna}):

\begin{defn}
A {\bf singular link} of degree m (within a linking class) is a link with m ordered pairs of 
double points, with both crossings in each pair involving the same 2 components.  The ordering 
of each pair is denoted by labeling the first crossing with a + and the second with a -.  I.e., 
each pair is of the form:
$$\psfig{file=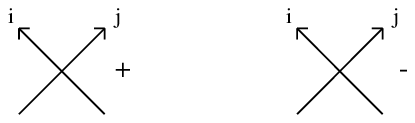}$$
\end{defn}

Given an invariant $V$ well defined within a linking class, we can extend it to singular links 
within that class by the relation:
$$\psfig{file=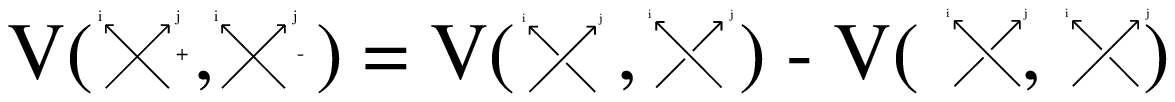}$$

\begin{defn}
An invariant $V$ defined within a linking class is of {\bf type m} if it is trivial on all 
singular links (within that linking class) of degree $\geq$ m+1.  $V$ is said to be of 
{\bf finite type} if it is of type m for some finite m.
\end{defn}

As with the clasping operations defined in the previous sections, we want to know what the 
equivalence classes of the double crossing change operation are.  The answer is given by the 
following proposition:

\begin{prop} \label{P:double-crossing-class}
If $L_1$ and $L_2$ are links with the same linking matrix, then $L_1$ can be transformed 
into $L_2$ (up to {\bf isotopy}) via double-crossing changes.
\end{prop}
{\sc Proof:}  $L_1$ and $L_2$ are the closures of string links $\sigma_1$ and $\sigma_2$.  
Clearly, it will suffice to prove the proposition for the associated string links.  Define 
the string link $\gamma = \sigma_2(\sigma_1)^{-1}$, so $\sigma_2 = \gamma\sigma_1$.  Since 
$\sigma_1$ and $\sigma_2$ have the same linking matrix, $\gamma$ is algebraically unlinked.  
It now suffices to show that $\gamma$ can be transformed into the trivial string link by 
double-crossing changes.

First consider components 1 and 2 of $\gamma$ (components are numbered from left to right).  
There are four ways in which these two components can cross (modulo a local rotation):
$$\psfig{file=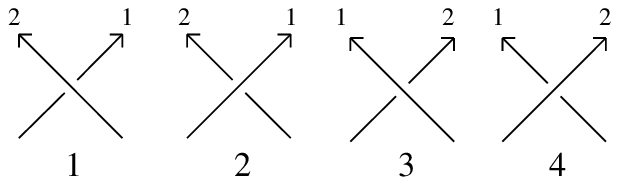}$$
These crossings appear with multiplicities $n_r, p_r, n_l, p_l$ respectively ($n, p$ refers 
to the sign of the crossing, $r, l$ refers to whether component 1 is moving to the right or 
the left).  Since the components are unlinked, $p_r+p_l = n_r+n_l$.  Also, since component 1 
must start and end to the left of component 2, $n_r+p_r = n_l+p_l$.  By taking the difference 
of these two equations, we find that $p_l = n_r$ and $p_r = n_l$.  So by changing the $p_r$ 
pairs of oppositely signed crossings (crossings of types 2 and 3 are paired), we are left with 
component 1 always undercrossing component 2.

By a similar argument for each pair of components, we can transform $\gamma$ so that component 
$i$ undercrosses component $j$ whenever $i<j$.  Hence, each component is at a different level, 
and the string link is now trivial.  $\Box$

\subsection{Chord Diagrams} \label{SS:double-crossing-chords}

We define chord diagrams as in the usual theory of finite type invariants, with the added 
condition that the chords come in ordered pairs.  Following \cite{bna}, we will call these 
diagrams {\it Double Dating Diagrams}:

\begin{defn}
A {\bf Double Dating Diagram (DD)} of degree m is a collection of l ordered oriented circles and 
m ordered pairs of lines (chords) so that both chords in each pair connect the same 
two circles.  The ordering of the pairs is denoted by labeling the first chord with a + and 
the second with a -, as seen below:
$$\psfig{file=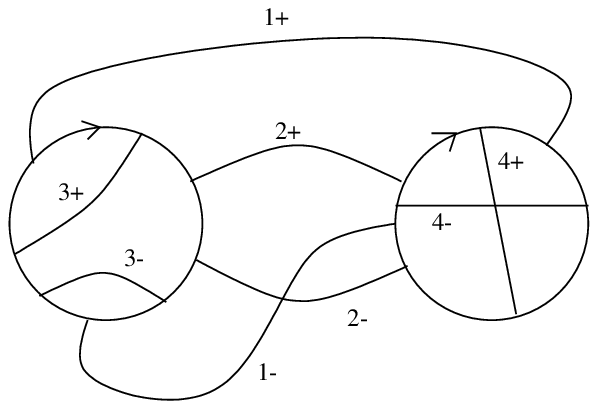}$$
\end{defn}

Given a singular link $L$ of degree $m$ (i.e. an immersion $L:\sqcup S^1 \rightarrow S^3$), 
there is a natural associated DD $D_L$ of degree $m$, 
where the pairs of chords of $D_L$ connect the preimages of the double points of $L$.  Conversely, 
given a DD $D$ of degree $m$ and a linking class, we can associate to $D$ a singular link $L_D$ 
in that linking class of degree $m$ by immersing $D$ in 3-space so that the two points joined 
by each chord are mapped to a double point of the link.  This link is not unique, but an 
argument similar to the proof of Proposition~\ref{P:double-crossing-class} shows that any two 
choices for $L_D$ (within a linking class) differ only by double-crossing changes (see also 
Theorem 2.1 of \cite{bna}).

This means that, given a finite type invariant $V$ of type $m$ in a given linking class, we can 
define a linear functional $W(V)$ on the space of DD diagrams of degree $m$ by the equation:
$$W(V)(D) = V(L_D)$$
This is well defined, because any two choices of $L_D$ differ by singular links of degree 
$\geq m+1$, on which $V$ is trivial.

We would now like to know some of the relations which $W(V)$ will satisfy:

\begin{prop} \label{P:chord-relations}
W(V) will satisfy the following relations:
\begin{itemize}
	\item  Antisymmetry relation:  Changing the order of a pair of chords (i.e. swapping 
the labels + and -) changes the sign of W(V).
		$$\psfig{file=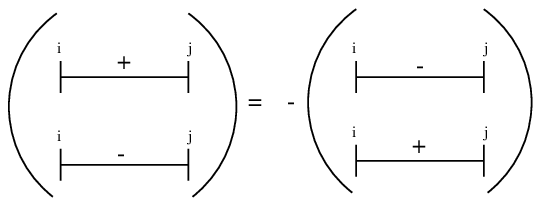}$$
	\item  Associative relation:
		$$\psfig{file=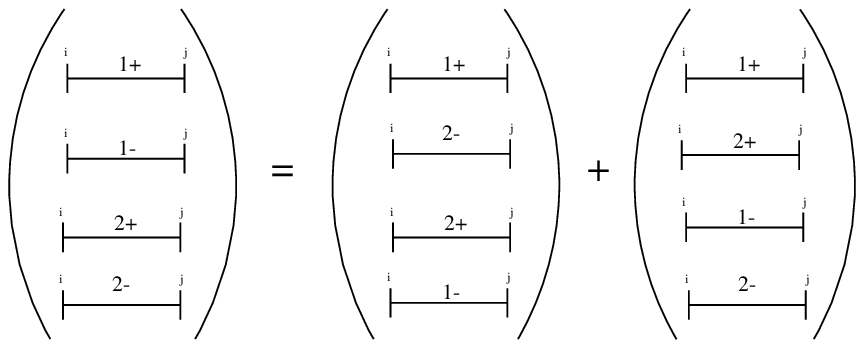}$$
	\item  1-term relation:
		$$\psfig{file=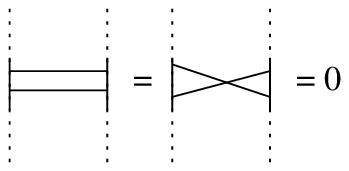}$$
	\item  4-term(a) relation:
		$$\psfig{file=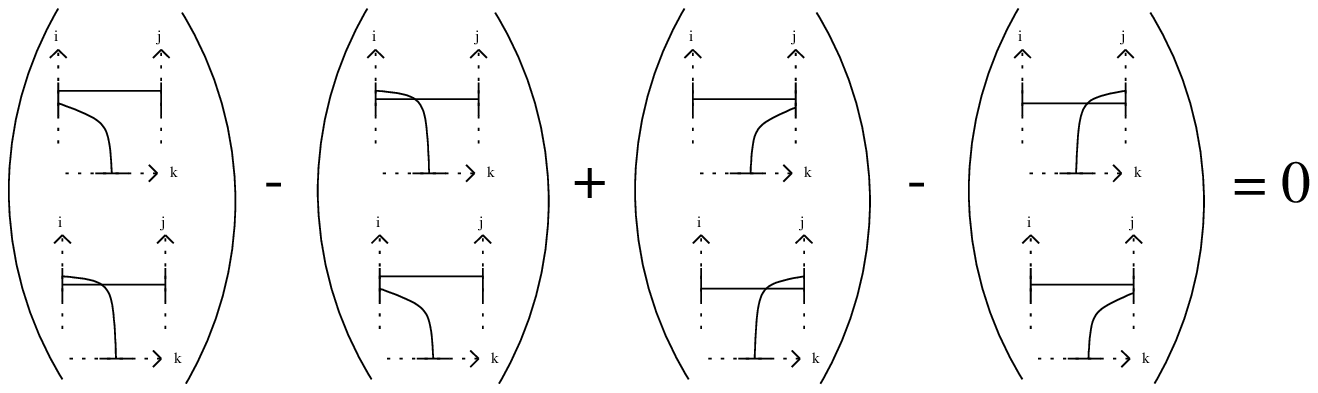}$$
	\item  4-term(b) relation:
		$$\psfig{file=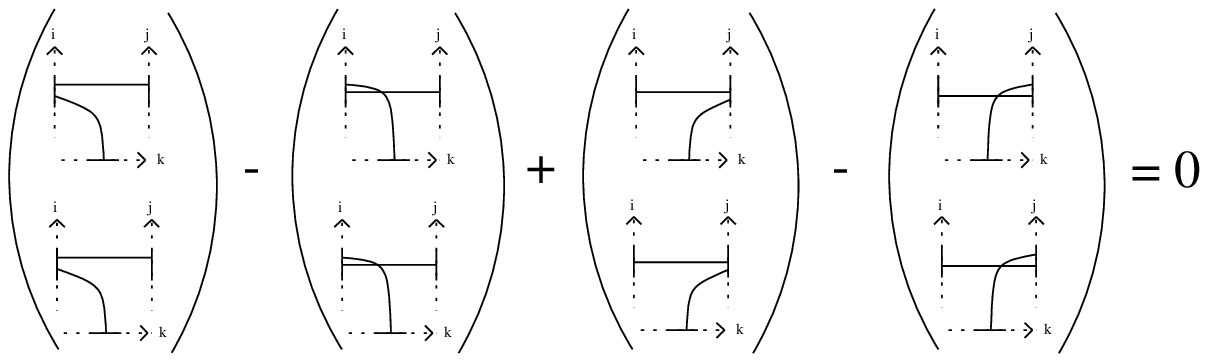}$$
\end{itemize}
By the Antisymmetry Relation, the labels (+/-) of the pairs in the last 3 relations are 
irrelevant.
\end{prop}
{\sc Proof:}  The Antisymmetry and Associative Relations are trivial (just expand the 
associated singular links).  The 1-term relation is a consequence of the second Reidemeister 
move:
$$\psfig{file=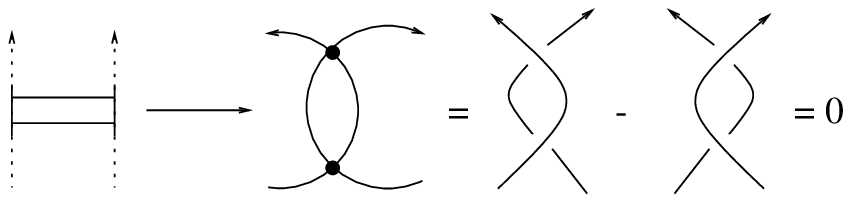}$$
The two 4-term relations are proved like the usual 4-term relation for chord diagrams 
(see \cite{bn1}), except that we are looking at two ``triple crossings'' at once.  We can 
view them as the result of bringing a contractible loop from underneath a double point to 
above it:
$$\psfig{file=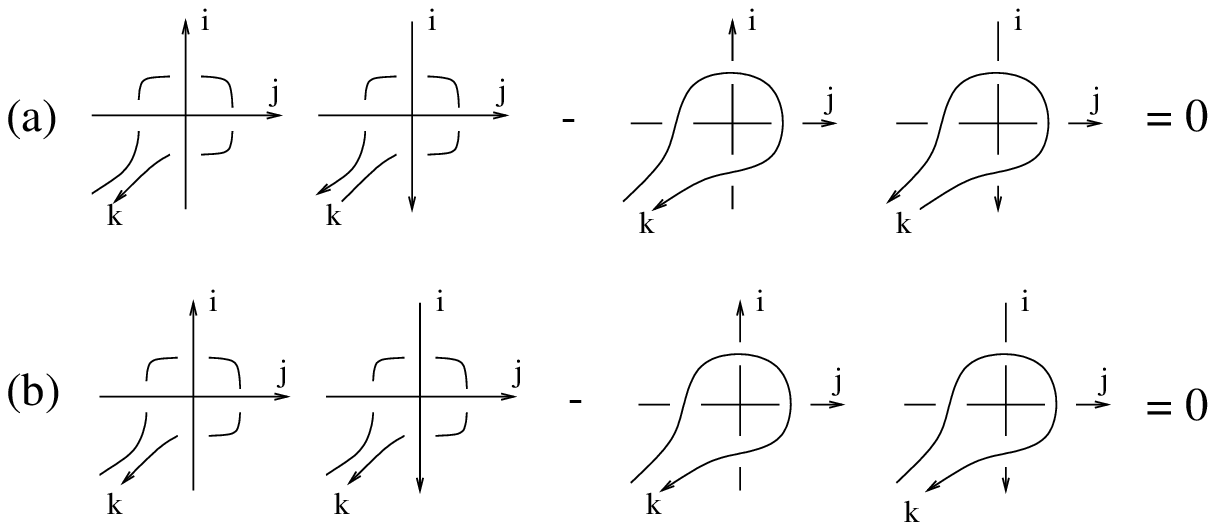}\Box$$

\subsection{$\bar{\mu}(ijk)$ is finite type} \label{SS:ijk-double-crossing}

In the next sections we will show that $\bar{\mu}(ijk)$ is, in each linking class, a finite type 
invariant of type 2; i.e. that it is trivial on singular links with 3 or more pairs of double 
points.  We will begin by showing that it is trivial on singular links with more than one 
pair of double points between the same two components, and use this result to prove the rest 
of the theorem.

We want to look at how $\bar{\mu}(ijk)$ is altered by double-crossing changes.  More generally, 
we want to see how crossing changes affect $\mu(ijk)$, before we mod out by $\Delta$.  Recall 
from \cite{mi} that $\mu(ijk)$ is invariant under a cyclic permutation of the indices, and 
changes sign if two indices are transposed.  So $\mu(ijk) = \mu(kij) = -\mu(jik)$.

Every time $l_k$ passes underneath another component of the link, it is multiplied by a generator 
of the Wirtinger presentation, which can be written as a word in the meridians of the link.  In 
fact, this word is just a conjugate of one of the meridians.  Consider the following picture:
$$\psfig{file=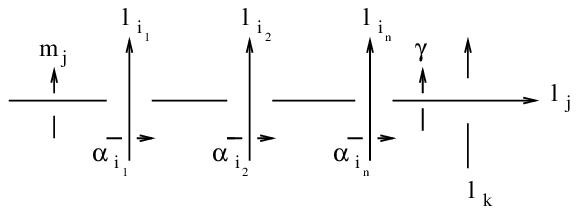}$$
The crossing of $l_k$ under $l_j$ inserts $\gamma$ into the word representing $l_k$.  $\gamma$ 
can be written as:
$$\gamma = \alpha_{i_n}^{-1}...\alpha_{i_1}^{-1}m_j\alpha_{i_1}\alpha_{i_n}$$
where $\alpha_{i_r} = \epsilon_r m_{i_r} \epsilon_r^{-1}$ for some word $\epsilon_r$.

The contribution to $\mu(ijk)$ by this crossing is (where $mult_i(w)$ is the multiplicity, with 
sign, of $m_i$ in the word $w$):
$$mult_i(\alpha_{i_n}^{-1}...\alpha_{i_1}^{-1})mult_j(m_j) 
+ mult_i(m_j)mult_j(\alpha_{i_1}\alpha_{i_n})$$
$$ = \left({\sum_{r = 1}^n{-mult_i(\alpha_{i_r})}}\right)mult_j(m_j) +
mult_i(m_j)\left({\sum_{r = 1}^n{mult_j(\alpha_{i_r})}}\right)$$
Since $mult_i(\alpha_{i_r}) = mult_i(\epsilon_r) + mult_i(m_{i_r}) - mult_i(\epsilon_r) 
= \delta(i,i_r)$, we find that the contribution is:
$$\left({\sum_{r = 1}^n{\delta(i,i_r)}}\right)mult_j(m_j) + 
mult_i(m_j)\left({\sum_{r = 1}^n{\delta(j,i_r)}}\right)$$
So it depends only on which components $l_j$ has previously passed under, not on {\it where} it 
passed under them, or in which order.

\begin{prop} \label{P:two-components}
If $L$ is a singular link with two pairs of double point between the same two components, then 
$\bar{\mu}(ijk)(L) = 0$.
\end{prop}
{\sc Proof:}  Let us denote the two pairs in question by $a = (a_+,a_-)$ and $b = (b_+,b_-)$.  
Then $L = L_{++} - L_{+-} - L_{-+} + L_{--}$, where the first index reflects the resolution of 
$a$, and the second index reflects the resolution of $b$.  Clearly, the only case of interest 
is when $a$ and $b$ are pairs of double points between two of the components $i,j,k$.  WLOG, 
we assume they are between $i$ and $j$.  Now we look at the crossings where $l_k$ passes 
under another component - these will be the same in all four links.  From the formulae above, 
we see that these crossings contribute to $\mu(ijk)$ only if $l_k$ is passing under $l_i$ or 
$l_j$.  The difference in the contributions of each crossing between $L_{++}$ and $L_{+-}$ is 
0, 1 or -1, depending on which of the double points $b_i$ occur between the basepoint of the 
component $l_k$ is undercrossing and the crossing in question.  The difference between 
$L_{-+}$ and $L_{--}$ will be the same, which means that $\bar{\mu}(ijk)(L)$ vanishes.  $\Box$

An immediate consequence of this proposition is the following corollary, which has also been 
proven (in a different way) by Appleboim and Bar-Natan in \cite{bna}.

\begin{cor} \label{C:ijk-finite-type}
$\bar{\mu}(ijk)$ is a finite type link homotopy invariant (within linking classes), of type 3.
\end{cor}
{\sc Proof:}  Since $\bar{\mu}(ijk)$ is a homotopy invariant, it is trivial on any singular 
link with a pair of double points between a component and itself.  So we only need to consider 
singular links where the double points are between {\it distinct} components of the link.  
Then $\bar{\mu}(ijk) = 0$ on any singular link with 4 or more pairs of double points, 
since either one of the pairs must involve a component other than $i,j,k$, or there will 
be two pairs connecting the same two components.  So $\bar{\mu}(ijk)$ is of finite type, of 
type 3.  $\Box$

\subsection{$\bar{\mu}(ijk)$ is type 2}

The only question which remains is whether $\bar{\mu}(ijk)$ is of type 1 or 2.  It cannot 
be of type 1 because of the following example:
$$\psfig{file=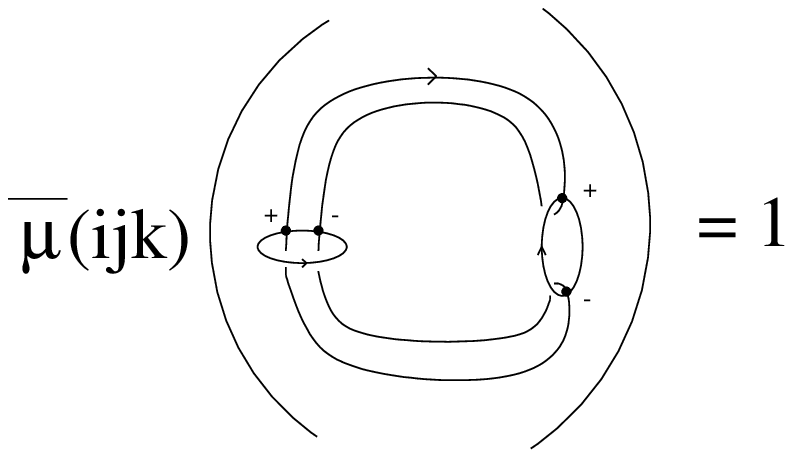}$$

We will show that $\bar{\mu}(ijk)$ is of type 2 by showing that the functional $W$ it induces on 
DD diagrams of degree 3 is trivial ($W$ is well-defined by Corollary~\ref{C:ijk-finite-type}).  
Since $\bar{\mu}(ijk)$ is a link homotopy invariant, $W$ 
will vanish on any DD diagram which has chords with both endpoints on the same component (see 
\cite{bn2}).  It is convenient to look at the DD diagrams for string links rather than links 
(just cut each circle to make an interval), remembering that chords are allowed to cycle from 
the top to the bottom.  Since $\bar{\mu}(ijk)$ only depends on components $i,j,k$, we only 
need to consider diagrams on 3 components.

So our first task is to count all the diagrams on 3 components with 3 pairs of chords.  By 
Proposition~\ref{P:two-components}, we can ignore diagrams with more than one pair of chords 
between the same 2 components (in these cases, we know $\bar{\mu}(ijk) = 0$).  So we will have 
one pair of chords connecting each of the 3 possible pairs of components.  This means there 
will be 4 endpoints on each component, allowing 4! = 24 permutations on each component.  Since 
we can cycle these without changing the diagram (since we are really looking at links), there 
are effectively 6 different permutations on each component, giving $6^3 = 216$ diagrams.  Since 
we want to show that $\bar{\mu}(ijk) = 0$ on each diagram, the signs of the diagram are 
irrelevant, so by the antisymmetry relation we can interchange the two chords in each pair 
with impunity.  This leaves us with $\frac{6^3}{2^3} = 3^3 = 27$ diagrams.  Actually, there 
are 28, because the last two are sufficiently symmetrical that permuting the endpoints of the 
chords cyclically yields only 32 different diagrams, rather than 64.  The 28 possible diagrams 
are listed in Figure~1.
	\begin{figure}
	$$\psfig{file=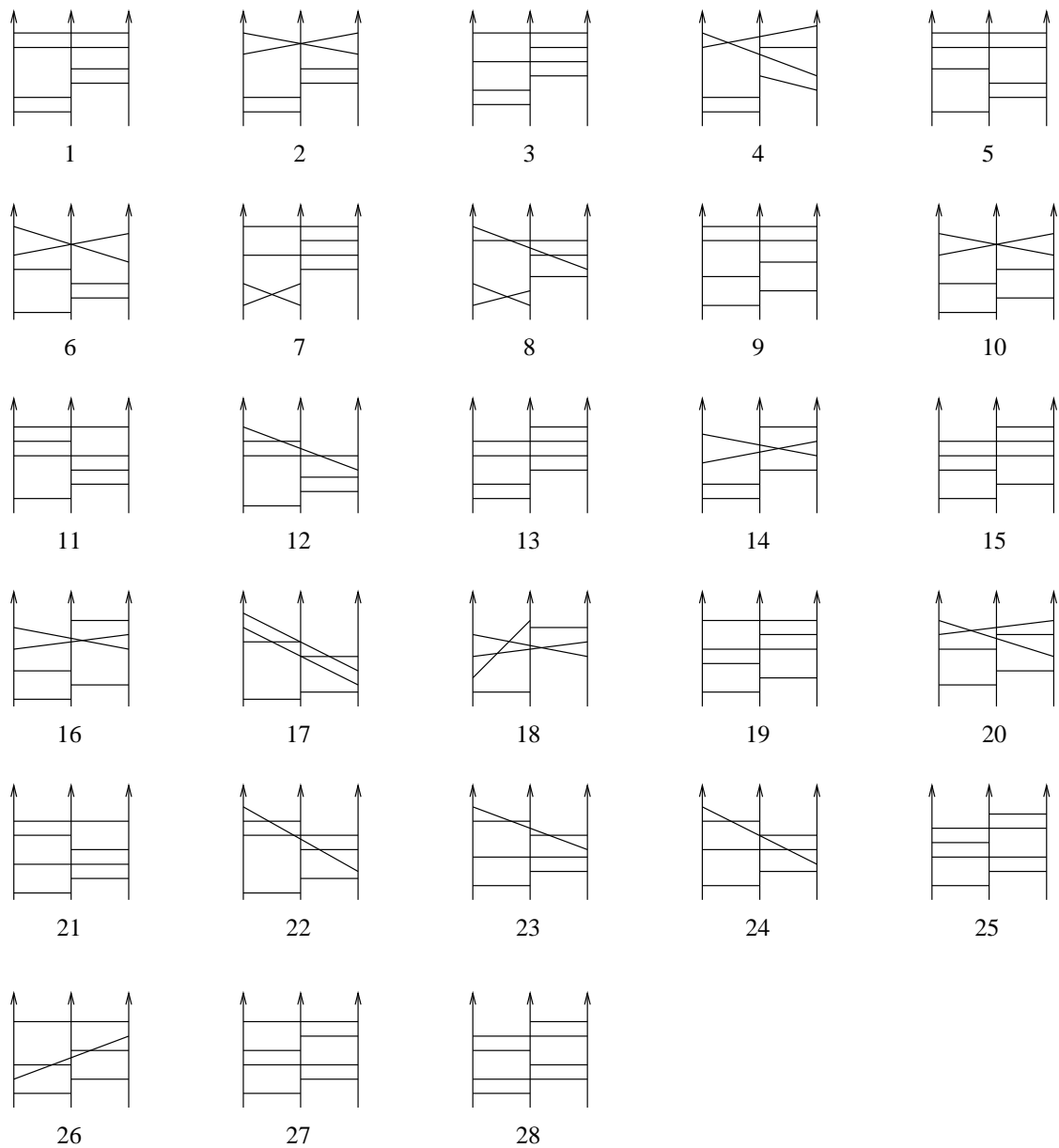}$$
	\caption{Diagrams of degree 3}
	\end{figure}
The first 18 diagrams are trivial by the 1-term relation.  By combining Proposition 3 with 
the 4-term and 1-term relations, we find that (as far as $\bar{\mu}(ijk)$ is concerned) the 
remaining 10 diagrams can be paired:  19 = 20, 21 = 22, 23 = 24, 25 = 26, 27 = 28.  For 
example, we will show how 19 = 20:
$$\psfig{file=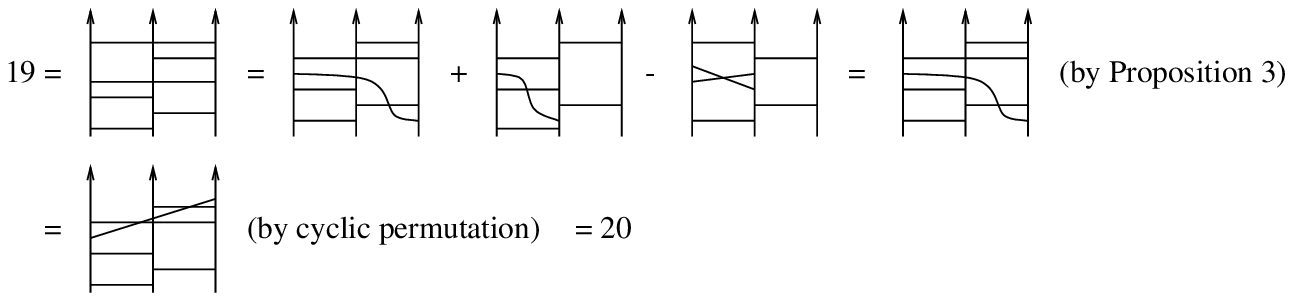}$$
Thus we are left with only 5 diagrams to consider.  By direct computation, we find that 
singular links representing these 5 diagrams in the trivial linking class (the class of 
algebraically split links) have $\bar{\mu}(ijk) = 0$.  This shows:

\begin{lem} \label{L:ijk-split-type2}
$\bar{\mu}(ijk)$ is of type 2 in the class of algebraically split links.
\end{lem}

We will use this fact to show that $\bar{\mu}(ijk)$ is of type 2 in every linking class.  

\begin{thm} \label{T:ijk-type2}
$\bar{\mu}(ijk)$ is a finite type invariant of type 2 in every linking class.
\end{thm}
{\sc Proof:}  Given a linking class, pick a link $S$ in that class, and pick a string link 
$\sigma_S$ such that $S = \hat{\sigma_S}$.  Then given a singular link $L$ of degree 3 in the 
class, and a singular string link $\sigma$ such that $L = \hat{\sigma}$, we can write 
$\sigma = \gamma\sigma_S$, where $\gamma$ is a singular link in the trivial linking class 
(i.e. algebraically unlinked) of degree 3.  Then $\gamma$ resolves into an alternating sum 
of 8 algebraically unlinked string links $\gamma_r$, and $L$ resolves into an alternating 
sum of 8 links $L_r$, where $L_r = \hat{\sigma_r} = \hat{\gamma_r\sigma_S}$.  By 
Lemmas~\ref{L:closing} and \ref{L:additivity}, $\bar{\mu}(ijk)(L_r) = \mu(\sigma_r)\ mod\ \Delta 
= \mu(ijk)(\gamma_r) + \mu(ijk)(\sigma_S)\ mod\ \Delta = \bar{\mu}(ijk)(\hat{\gamma_r}) + 
\bar{\mu}(ijk)(S)\ mod\ \Delta$.  Then:
$$\bar{\mu}(ijk)(L) = \sum_{r=1}^8{(-1)^r(\bar{\mu}(ijk)(\hat{\gamma_r})+\bar{\mu}(ijk)(S))}
\ mod\ \Delta$$
The second terms all cancel, so $\bar{\mu}(ijk)(L) = \bar{\mu}(ijk)(\hat{\gamma})$.  Since 
$\bar{\mu}(ijk)$ is of type 2 in the class of algebraically split links, this is 0.  Hence, 
$\bar{\mu}(ijk)$ is of type 2 in every linking class.  $\Box$

\section{Comparing Theories of Finite Type Invariants} \label{S:compare}

In this section we want to compare our two theories of finite type invariants in linking classes - 
one generated by Borromean clasps, and the other by double crossing changes.  First we will make 
this comparison precise by introducing the idea of {\it local equivalence} of operations (also 
see \cite{na}).

We will consider two operations $A$ and $B$ on links.  We will assume that these 
are both {\it local} moves in the sense that each takes place within a small ball (or a finite 
collection of small balls), leaving the rest of the knot fixed.

\begin{defn}
A {\bf locally generates} B if any move B, in some small neighbourhood, is the result of a finite 
sequence of moves A in the same neighbourhood.  In particular, since this is a purely local 
criterion, the number of moves A required is fixed.
\end{defn}

\begin{defn}
A and B are {\bf locally equivalent} if each locally generates the other.
\end{defn}

Any local move will generate a theory of finite type invariants in the obvious way (as we have 
done earlier in this paper).  We will call the theory generated by a local move $A$ the $A$-
theory.

\begin{prop} \label{P:locally-generates}
If A locally generates B, and V is a finite type link invariant of type n in the A-theory, 
then V is of type n in the B-theory.
\end{prop}
{\sc Proof:}  Consider a link $L$ with $n+1$ $B$-singularities.  Since each $B$-move is the 
result of a fixed number (say $k$) $A$-moves, each $B$-singularity is (locally) a linear 
combination of $k$ $A$-singularities.  Hence $L$ equals a linear combination of $k^{n+1}$ links, 
each with $n+1$ $A$-singularities.  So $V(L) = V(sum) = 0$, so $V$ is also of type $n$ in the 
$B$-theory.  $\Box$

\begin{cor} \label{C:locally-equivalent}
If A and B are locally equivalent, any invariant of type n in the A-theory is also of type n 
in the B-theory.  We will say that the theories are {\bf isomorphic}.
\end{cor}

Now we can consider the particular examples of the Borromean clasp theory and the double 
crossing change theory.  The first difference to note is that we have shown that any two 
links in the same linking class are equivalent up to isotopy modulo double crossing changes, 
but only up to homotopy modulo Borromean clasps.  In fact, this is not a problem - we can 
strengthen our result for Borromean clasps (though not, I think, for the higher-order 
clasps).

\begin{prop} \label{P:borr-class-isotopy}
Any two links with the same linking numbers are equivalent (up to {\bf isotopy}) modulo 
Borromean clasps.
\end{prop}
{\sc Proof:}  Murakami and Nakanishi note in \cite{mn} that adding a Borromean clasp locally 
generates their $\Delta$ unknotting operation, shown below (for oriented links, we consider 
this operation with all possible orientations, giving 8 oriented operations):
$$\psfig{file=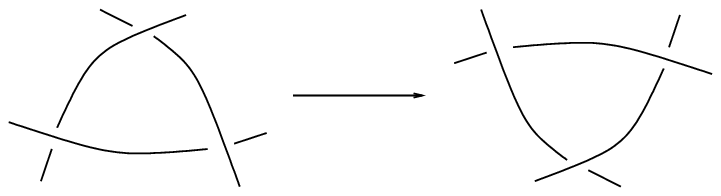}$$
It is easy to see that the $\Delta$ move also locally generates the Borromean clasp, so the 
two moves are locally equivalent.  Murakami and Nakanishi show that any two links in the 
same linking class (in their terminology, link-homologous) are equivalent, up to isotopy,  
modulo $\Delta$ moves (in particular, any component can be unknotted).  Hence, they are also 
equivalent up to isotopy modulo Borromean clasps.  $\Box$

Now we can state the main result of this section:

\begin{thm} \label{T:clasps contain double-crossings}
If V is a finite type invariant in the double crossing change theory of type n, then V is also 
of finite type in the Borromean clasp theory, and also of type n (though possibly also of lower 
type).
\end{thm}
{\sc Proof:}  It is obvious that Borromean clasps are locally generated by double crossing changes.  
The result is then given by Proposition~\ref{P:locally-generates}.  $\Box$

However, the converse of this theorem is false.  $\bar{\mu}(ijk)$ is of type 1 in the Borromean 
clasp theory, but it is {\it not} of type 1 in the double crossing change theory.  So it appears 
that there may be more finite type invariants in the Borromean clasp theory.

{\sc Remark:}  Since the usual crossing change operation locally generates both double crossing 
changes and Borromean clasps, any of the usual finite type invariants are also of finite type 
in these theories.  This has also been shown (for the double crossing change theory) in 
\cite{bna}.

\section{Questions} \label{S:questions}

In section~\ref{S:claspclass} we described the equivalence classes of the clasping operations at 
all levels.  Of course, we would really like these classes to be the $n$-linking classes.  
Certainly, if two string links are 
equivalent modulo $H_n(k)$, then they will be indistinguishable by $\mu$-invariants of 
length $n$ or less, and so their closures will be in the same $n$-linking class.  However, the 
converse is not so clear, and leads to the following question:

\begin{quest} \label{Q:claspclass=linkingclass}
Say that L and $L'$ are two k-component links in the same $n$-linking class (so that 
$\bar{\mu}(i_1...i_r)(L) = \bar{\mu}(i_1...i_r)(L')$ for every $\{i_1,...,i_r\}$ with $r\leq n$).  
Do there exist string links $\sigma, \sigma'$ such that L = $\hat{\sigma}$, $L' = \hat{\sigma}'$, 
and $\sigma \equiv \sigma'$ modulo $H_n(k)$?
\end{quest}

Essentially, this is asking whether the indeterminacy $\Delta$ of the $\bar{\mu}$-invariants is 
the ``same'' as the subgroup of the group of conjugations and partial conjugations which preserve 
the closure of the string link $\sigma$ and also fix it modulo $H_n(k)$.

Several questions arise from the comparison of the Borromean clasp and double crossing change 
theories in section~\ref{S:compare}.  We have shown that the double crossing change theory is in some 
sense a ``subset'' by the Borromean clasp theory, but we don't know if it is a {\it proper} 
subset:

\begin{quest} \label{Q:proper-subset}
Is there a finite type invariant in the Borromean clasp theory which is not of finite type 
in the double crossing change theory?
\end{quest}

On a broader level, we would like to know which is the ``better'' theory of finite type invariants 
in linking classes.  So far, it seems that the Borromean clasp theory may be better, as 
potentially providing more invariants to work with.  However, much of the work on finite type 
invariants has been done by studying chord diagrams, and it is not clear what these would be 
in the Borromean clasp theory:

\begin{quest} \label{Q:borr-chord-diagrams}
Is there a useful graded vector space of ``chord diagrams'' associated to singular links in the 
Borromean clasp theory such that any two singular links associated to the same chord diagram 
are equivalent modulo Borromean clasps?  In other words, does any link invariant give rise to a 
well-defined weight system?
\end{quest}

The problem here is whether the Borromean clasps would need to reach ``inside'' the singularities, 
which would not be allowed.

Finally, we can ask about finite type {\it concordance} invariants in linking classes.  As we 
have noted, the Borromean clasp and double crossing change theories both apply equally well 
to considering links (within linking classes) up to isotopy (and hence concordance), rather 
than homotopy.  So we can also ask whether Milnor's {\it concordance} invariants (when the 
indices repeat) are of finite type.  

\begin{quest}
Is $\bar{\mu}(i_1...i_n)$ of finite type (in some sense), when the $i_j$'s can repeat?
\end{quest}

\section{Acknowledgements}

I want to particularly acknowledge the help of Paul Melvin, who gave me the idea of working 
within linking classes and helped me begin.  I also want to thank Eli Appleboim for allowing 
me to see an early draft of his work.  Finally, I want to thank Rob Kirby, and all the 
members of Berkeley's Informal Topology Seminar, for their advice and criticisms.

\end{document}